\documentclass[12pt]{article}
\usepackage{latexsym,bm,amssymb,amsfonts,amsbsy,amsmath,graphics,graphicx,color}
\usepackage{url}
\newtheorem{theo}{Theorem}

\def\RR{\mathbb R}

\def\pmatrix{ \left( \begin{array} }
\def\endpmatrix{ \end{array} \right) }

\def\dd{\mathrm{d}}

\def\II{{\cal I}}
\def\P{{\cal P}}

\def\d2dxx{\frac{\partial^2}{\partial x^2}}

\def\no{\noindent}
\def\diag{{\rm diag}}

\def\phi{\varphi}

\def\P{{\cal P}}

\begin{document}

\title{Energy-conserving methods for Hamiltonian Boundary Value Problems and applications in astrodynamics.}

\author{Pierluigi Amodio$^a$        \and
        Luigi Brugnano$^b$  \and
        Felice Iavernaro$^a$
}

\date{\small$^a$ Dipartimento di Matematica\\ Universit\`a di Bari, Italy\\[2mm]
$^b$ Dipartimento di Matematica e Informatica ``U.\,Dini''\\ Universit\`a di Firenze, Italy}

\maketitle

\begin{abstract}
We introduce new methods for the numerical solution of general Hamiltonian boundary value problems. The main feature of the new formulae is to produce numerical solutions along which the energy is precisely conserved, as is the case with the analytical solution. We apply the methods to locate periodic orbits in the circular restricted three body problem by using their energy value rather than their period as input data. We also use the methods for solving optimal transfer problems in astrodynamics. 

\bigskip
\no {\bf Keywords:} Energy conserving  Runge-Kutta methods, Hamiltonian boundary value problems, Astrodynamics, Optimal control.

\bigskip
\no{\bf PACS:} 02.30.Hq, 02.60.Lj, 45.20.dh, 45.20.Jj.

\bigskip
\no{\bf MSC:} 65P10, 65L10 65L06.
\end{abstract}

\section{Introduction}
\label{intro}
We are concerned with the numerical solution of the general autonomous Hamiltonian boundary value problem
\begin{equation} 
\label{ham_bvp}
\left\{ 
\begin{array}{ll}
\dot y(t) = J\nabla H(y(t)), \quad t\in[t_0,t_f] \\ 
g(y(t_0),y(t_f))=0.
\end{array}
\right.
\end{equation}
The scalar function $H:\Omega\subset \RR^{2m} \rightarrow \RR$ is the {\em Hamiltonian} of the problem, $J=\pmatrix{cc} 0 &I_m\\ -I_m&0\endpmatrix$ (here and in the sequel $I_r$ will denote  the identity matrix of dimension $r$), and the vector function $g:\RR^{2m}\times \RR^{2m} \rightarrow \RR^{2m}$ defines the boundary conditions. Hereafter, both $H$ and $g$ will be assumed to be suitably regular. As is well known, the  value of the Hamiltonian function is constant along the solution of \eqref{ham_bvp}. An easy manner to see this is to consider the line integral associated with the vector field $\nabla H(y)$ evaluated along the path defined by the solution $y(t)$ of  \eqref{ham_bvp}, which equals the variation of $H$ along the end-points of the path. Exploiting the skew-symmetry of matrix $J$, we have, for $t_0\le t \le t_f$,
\begin{equation}
\label{lin_int}
\begin{array}{rl}
\displaystyle H(y(t))-H(y(t_0)) & \displaystyle = \int_{t_0}^{t} \dot y^T(\tau) \nabla H(y(\tau)) \mathrm{d} \tau \\
&\displaystyle  = \int_{t_0}^{t} \nabla^T H(y(\tau)) J^T  \nabla H(y(\tau)) \mathrm{d} \tau =0.
\end{array}
\end{equation} 
The state vector $y$ splits in two vectors of length $m$, $y^T=[q^T,p^T]$ referred to as generalized coordinates and conjugate momenta. The  numerical treatment of Hamiltonian problems is thoroughly discussed in the monographs \cite{HLW,LR,SC}.

The aim of the present work is to construct energy-conserving methods for problem \eqref{ham_bvp}, that is methods producing numerical solutions $\{y_i\}$ along which the value of the energy is the same: $H(y_i)=H(y_{i-1})$. 

Interest in problems such as \eqref{ham_bvp} arises in several research areas. In this paper (see Section~\ref{sec:4}), we focus our attention on some applications in celestial mechanics and astrodynamics. In particular, we consider the dynamics of a massless object (planetoid) subject to the gravitational field induced by two massive bodies (primaries) revolving in circular orbits about their center of mass. Such a dynamical system, referred to as the \textit{circular restricted three-body problem}, together with its generalizations, has been deeply studied since Poincar\'e. Its renewed interest is motivated by the fundamental role it plays in the context of space mission design and control problems in aerospace engineering, such as the nonlinear trajectory optimization  and the spacecraft orbit transfer \cite{Ba,KLMR}. 

 The paper is organized as follows. In the next section we briefly recall the definition of the energy-conserving methods named HBVMs and describe their main features. For a detailed description of HBVMs, and their properties when applied to Hamiltonian IVPs, see \cite{BIT2,BIT3,BIT4}.  The implementation of HBVMs to solve problem \eqref{ham_bvp} will be discussed in Section~\ref{sec:3}. Section~\ref{sec:4} will be devoted to the description of specific problems related to the circular restricted three-body system, and their numerical treatment. A few concluding remarks are then reported in Section~\ref{sec:5}. For an introduction on the solution of general  boundary value problems by using one-step methods see, for example, \cite{AMR}.    


\section{Definition of the methods}
\label{sec:2}
In this section we recall the definition of HBVMs. These are Runge-Kutta methods characterized by a low-rank coefficient matrix. 

The main prerogative of a HBVM is to reproduce, in the discrete setting, property \eqref{lin_int} of conservative vector fields. To this end, we consider the approach discussed in \cite{BIT4} and exploiting a Fourier expansion of the continuous problem $\dot y(t) =J\nabla H(y(t))$ restricted to the interval $t \in [t_0,t_0+h]$, where $h=\frac{t_f-t_0}{n}>0$ will act as the stepsize of integration in the one-step method that will finally arise from this analysis. The procedure is then iterated on adjacent intervals $[t_i,t_{i+1}]$ with $t_{i+1}=t_i+h$, $i = 0,\dots,n-1$, until the overall integration interval $[t_0,t_f]$ is covered. 

Let us then consider the {\em Legendre} polynomials $P_i$ shifted on the interval $[0,1]$, and scaled in order to be orthonormal:
\begin{equation}\label{orto}
\deg P_i = i, \qquad \int_0^1 P_i(x)P_j(x)\dd x = \delta_{ij}, \qquad \forall i,j\ge 0,
\end{equation}
where $\delta_{ij}$ is the Kronecker symbol. The roots $\{c_1,\dots,c_k\}$ of $P_k(x)$ are all distinct and symmetrically distributed on the interval $(0,1)$. Usually, they are referred to as the Gauss-Legendre abscissae on $[0,1]$ and generate the well-known Gauss-Legendre quadrature formulae, whose weights we denote $b_i$, $i=1,\dots,k$. The infinite sequence $\{P_i(t)\}$ forms an orthonormal basis of $L^2([0,1])$. Expanding the right-hand side of \eqref{ham_bvp} along this basis and truncating the series after $s$ terms changes the original differential equation $\dot y(t) =J\nabla H(y(t))$ to 
\begin{equation}
\label{fourier}
\dot \omega(t_0+ch) =  \sum_{j=0}^{s-1} P_j(c) \int_0^1 P_j(x) J\nabla H(\omega(xh)) \mathrm{d} x, \qquad c\in [0,1].
\end{equation} 
Notice that the solution $\omega(t_0+ch)$ of \eqref{fourier}  is indeed a polynomial of degree $s$. In \cite[Theorem 1]{BIT4} it has been shown that $$y(t_0+h)-\omega(t_0+h)=O(h^{2s+1}), \qquad \forall k\ge s,$$ so that, iterating the procedure \eqref{fourier} sequentially over the intervals $[t_i,t_{i+1}]$, $i=0,\dots,n-1$,  provides an approximation of order $2s$ to the true solution, on the whole interval $[t_0,t_f]$. One interesting aspect of formula \eqref{fourier} is that it inherits the energy conservation property of the original problem. In fact, by setting 
\begin{equation}
\label{gammaj}
\gamma_j(\omega)=\int_0^1 P_j(x) J\nabla H(\omega(xh))\mathrm{d}x, \quad j=0,\dots,s-1,
\end{equation} 
we have
\begin{equation}
\label{conservation_omega}
\begin{array}{rl}
\displaystyle H(\omega(t_1))-H(\omega(t_0)) & \displaystyle = h\int_{0}^{1} \dot \omega(t_0+ch)^T \nabla H(\omega(t_0+ch)) \mathrm{d} t \\
&\displaystyle  = \sum_{j=0}^{s-1} \gamma_j(\omega)^T J^T \gamma_j(\omega)=0.
\end{array}
\end{equation}
A drawback of formula \eqref{fourier} is the presence of the integrals defining the scalar products \eqref{gammaj}, which make it unusable for a direct implementation. To circumvent this problem, a quadrature formula is introduced to approximate such integrals. In particular, we consider the Gauss-Legendre quadrature based at $k\ge s$ abscissae $c_i$ and weights $b_i$, $i=1,\dots,k$, introduced earlier, thus obtaining the following approximation to \eqref{gammaj}:
\begin{equation}
\label{gammacj}
\hat \gamma_j(\omega)\equiv \sum_{\ell=1}^k b_\ell P_j(c_\ell) J\nabla H(\omega(c_\ell)) = \gamma_j(\omega)+O(h^{2k-j}).
\end{equation} 
In so doing, the polynomial approximation changes, so that we finally arrive at the method 
\begin{equation}\label{y1sig}
y_1=\Phi_h(y_0)\equiv \sigma(t_0+h),\end{equation} where the polynomial $\sigma\in\Pi_s$ is defined as (compare with \eqref{fourier})
\begin{equation}
\label{hbvm}
\dot \sigma(t_0+ch) =  \sum_{j=0}^{s-1} P_j(c) \sum_{\ell=1}^k b_\ell P_j(c_\ell) J\nabla H(\sigma(c_\ell)).
\end{equation}
It can be shown that this method has still order $2s$. However, due to the approximation \eqref{gammacj}, the polynomial $\sigma$ does not retain, in general, the conservation property \eqref{conservation_omega} of the original polynomial $\omega$. Neverthels, this is no much of an issue since the following two situations may occur:
\begin{enumerate}
\item $H(y)$ is a polynomial of degree, say $\nu$. In this case, the integrand in \eqref{gammaj} has degree at most $\nu s-1$ and, since the Gauss-Legendre quadrature formula is exact for polynomials of degree at most $2k-1$, it will be enough to choose $k\ge \frac{\nu s}{2}$ to get $\hat \gamma_j = \gamma_j$ and hence energy conservation. Indeed, in such a case one evidently obtains $\omega\equiv\sigma$;

\medskip
\item $H(y)$ is a general, though suitably regular, non-polynomial function. According to the analysis in \cite{BIT4}, one then proves that (see \eqref{y1sig}--\eqref{hbvm})
\begin{equation}\label{errH}H(y_1)-H(y_0) = O(h^{2k+1}).\end{equation}
Consequently, even in this case, we can get a \textit{practical energy conservation} by choosing $k$ as large as to guarantee that the error $O(h^{2k+1})$, appearing at the right-hand side in \eqref{errH}, is of the order of the machine epsilon. As we will see in the next section, choosing a large $k$ does not affect the overall computational cost associated with the implementation of the method, which essentially depends on $s$.
\end{enumerate}

Integrating both sides of \eqref{hbvm} with respect to the variable $c$ and evaluating at $c=c_i$ yields
\begin{equation}
\label{hbvm1}
\sigma(t_0+c_ih) = y_0 +h \sum_{j=0}^{s-1} \int_0^{c_i} P_j(\tau)\mathrm{d} \tau \sum_{\ell=1}^k b_\ell P_j(c_\ell) J\nabla H(\sigma(c_\ell)),\quad i=1,\dots,k.
\end{equation}
This formula is called \textit{Hamiltonian Boundary Value Method} (HBVM) and is tantamount to a Runge-Kutta collocation-like method with internal stages  $Y_i\equiv \sigma(t_0+c_ih)$. In fact, by setting 
$$\P_s =\pmatrix{ccc}
P_0(c_1) & \dots & P_{s-1}(c_1)\\
\vdots   &           &\vdots \\
P_0(c_k) &  \dots &P_{s-1}(c_k)\endpmatrix, \quad
\II_s = \pmatrix{ccc}
\int_0^{c_1}P_0(x)\dd x & \dots &\int_0^{c_1}P_{s-1}(x)\dd x\\
\vdots   &           &\vdots \\
\int_0^{c_k}P_0(x)\dd x &  \dots &\int_0^{c_k}P_{s-1}(x)\dd x\endpmatrix 
$$
$c=(c_1,\dots,c_k)^T$, $b=(b_1,\dots,b_k)^T$ and $\Omega = \diag(b_1, b_2,\dots, b_k)$, \eqref{y1sig}-\eqref{hbvm1} is equivalent to the $k$-stage R-K method defined by the following Butcher tableau:
\begin{equation}
\label{hbvm-rk}
\begin{array}{c|c}
c & \II_s\P_s^T\Omega \\
\hline 
& b^T
\end{array}.
\end{equation}
This method is denoted by HBVM$(k,s)$ to outline its dependence on the two integers $s$ (degree of the polynomial approximation, clearly determined by number of Legendre polynomials involved) and $k$ (which is related to the number of internal abscissae and, therefore, to the order of the quadrature). 

Notice that $\P_s$ and  $\II_s$ are ${k\times s}$ matrices while $\Omega\in R^{k\times k}$. 
When $k=s$ one can show that \eqref{hbvm-rk} becomes the usual Gauss collocation method of order $2s$ \cite{BIT2}. For any $k\ge s$ the coefficient matrix $A=\II_s\P_s^T\Omega$ has constant rank $s$ and its nonzero eigenvalues coincide with those of the Butcher matrix defining the basic $s$-stage Gauss method \cite{BIT3}. The main properties of HBVMs are summarized in the following theorem (see \cite[Corollary 3]{BIT4}). 
\begin{theo} HBVM$(k,s)$ is symmetric, of order $2s$ and energy-conserving for all polynomial Hamiltonians of degree $\nu\le \frac{2k}{s}$.
In any other case, $H(y_1)-H(y_0)=O(h^{2k+1})$, provided that $H$ is suitably regular.\end{theo}

\section{Simplified Newton iteration and implementation details}
\label{sec:3}

That the rank of the coefficient matrix $\II_s\P_s^T\Omega$ is $s$, independently of $k$, suggests that $k-s$ stages $Y_i$ may be regarded as linear combinations of the remaining $s$ stages. This algebraic property turns out to be of fundamental importance to reduce the computational effort associated with the implementation of the method when applied to \eqref{ham_bvp}.  Therefore, it is convenient to derive an alternative (though equivalent)  shape of an HBVM$(k,s)$ method.

Notice that the polynomial $\sigma$ in \eqref{hbvm} has degree $s$.  Hence, it is completely determined by $s$ (rather than $k$) stages, plus the condition $\sigma(t_0)=y_0$.  These $s$ stages have been called \textit{fundamental stages} and, without loss of generality, in order to maintain the notation as simple as possible, they are assumed to be the first ones: $Y_i$, $i=1,\dots,s$.\footnote{In the actual implementation, their distribution is chosen according to what explained in \cite{BIT1}, i.e., the corresponding $s$ abscissae are approximately uniformly spaced in $[0,1]$.} The remaining stages $Y_j$, $j=s+1,\dots,k$, though contributing in defining the final shape of the polynomial $\sigma$, may be conveniently defined as linear combinations of the fundamental stages, by simply setting $Y_j=\sigma(t_0+c_jh)$. In other words, setting $Z\equiv (Y_1^T,\dots,Y_s^T)^T$, $W\equiv (Y_{s+1}^T,\dots,Y_k^T)^T$, we have 
\begin{equation}\label{silsta}
W=a_0\otimes y_0+A\otimes I_{2m} Z, 
\end{equation}
where the entries of the vector $a_0\in\RR^{k-s}$ and the matrix $A\in\RR^{k-s\times s}$ are the evaluations, at the abscissae $c_{s+1},\dots,c_k$, of the Lagrange polynomials defined on the nodes $\{0, c_1,\dots,c_s\}$. For this reason, the stages in $W$ have been referred to as \textit{silent stages} \cite{IaTr09}. In conclusion, we arrive at the following formulation of the method
\begin{equation}
\label{hbvm2}\left\{
\begin{array}{rcl}
- e\otimes y_0 + Z &=& h(B_1\otimes J) \nabla H(Z) + h(B_2\otimes J) \nabla H(W),\\
W &=& a_0\otimes y_0+A\otimes I_{2m} Z, \\
-y_0+y_1 &=& h ( \beta_1^T \otimes J) \nabla H(Z)+h (\beta_2^T \otimes J) \nabla H(W).
\end{array}
\right.
\end{equation}
where the first block-equation corresponds to the first $s$ equations in \eqref{hbvm1}, the second block-equation, defining the silent stages, is inherited from (\ref{silsta}), and the last equation defines the new approximation, $y_1$, having set  $\beta_1=[b_1,\dots,b_s]^T$ and $\beta_2=[b_{s+1},\dots,b_k]^T$. The clear advantage of \eqref{hbvm2}, with respect to \eqref{hbvm-rk}, is that now the nonlinear and linear part of the system defining the stages are completely uncoupled.

Suppose that the interval $[t_0,~t_f]$ is divided into $n$ equispaced sub-intervals $[t_{i-1}, t_i]$, $i=1,\dots,n$, of length $h$. Then, equation \eqref{hbvm2} may be subsequently iterated on such intervals to yield the approximations $y_i\simeq y(t_i)$, 
for $i=0,\dots,n$.
In particular, $y_0, y_1, \dots, y_n$ are combined with the given boundary conditions to yield a large nonlinear system in the unknowns $y_0, Z_0, y_1, Z_1, \dots,$ $y_{n-1}, Z_{n-1}, y_n$, where $Z_i$ is the block-vector of the fundamental stages   (denoted by $Z$ in the first equation of \eqref{hbvm2}) associated with $y_i$.

Ignoring momentarily the boundary conditions, a Newton-like iteration applied to the nonlinear equations gives, for $i=0,\dots, n-1$, the sequences
$\{y_i^{(j)}, Z_i^{(j)} \}$ defined as (for sake of brevity, let us assume $Z_n^{(j)}\equiv0$)
$$
\left\{
\begin{array}{ll}
y_i^{(j+1)} &= y_i^{(j)}+\delta_i^{(j)}, \\
Z_i^{(j+1)} &= Z_i^{(j)}+\Delta_i^{(j)},
\end{array}
\right.$$
where the increments $\delta_i^{(j)}$ on the $y_i$ variable and $\Delta_i^{(j)}$  on the $Z_i$ variable are the solution of the following linear system with sparse structured coefficient matrix:
\begin{equation}
\label{struct_matr}
\pmatrix {ccccccccc}
V_1 & K_1 \\
L_1 & U_1^T & I_{2m} \\
&& V_2 & K_2 \\
&& L_2 & U_2^T & I_{2m} \\
&&&& \ddots & \ddots \\
&&&&&& V_n & K_n \\
&&&&&& L_n & U_n^T & I_{2m}
\endpmatrix
\pmatrix {c}
\delta_0^{(j)} \\ \Delta_0^{(j)} \\ \delta_1^{(j)} \\ \Delta_1^{(j)} \\ \delta_2^{(j)} \\ \vdots \\ \delta_{n-1}^{(j)} \\ \Delta_{n-1}^{(j)} \\ \delta_n^{(j)}
\endpmatrix =
\pmatrix {c}
b_1 \\ c_1 \\ b_2 \\ c_2 \\ \vdots \\  b_n \\ c_n
\endpmatrix.
\end{equation}
The blocks $V_i,U_i\in \RR^{2ms \times 2m}$, $L_i\in \RR^{2m \times 2m}$ and $K_i \in \RR^{2ms \times 2ms}$ are defined as follows:\footnote{In order not to complicate the notation, and for sake of brevity, we shall omit the iteration index $j$, for these blocks.}
\begin{equation}
\label{blocks}
\begin{array} {ll}
V_i &=-e \otimes I_{2m} - h(B_2 a_0) \otimes J \nabla^2 H(\bar y_i^{(j)}), \\
U_i^T &= -h (\beta_2^T A + \beta_1^T) \otimes J \nabla^2 H(\bar y_i^{(j)}), \\
L_i &= -I_{2m} - h (\beta_2^T a_0) J \nabla^2 H(\bar y_i^{(j)}), \\
K_i &= I_s \otimes I_{2m} - h (B_2 A + B_1) \otimes J \nabla^2 H(\bar y_i^{(j)}),
\end{array}
\end{equation}
where, for symmetry reasons,  $\bar y_i^{(j)}=\frac{y_i^{(j)}+y_{i-1}^{(j)}}{2}$. Finally, $b_i\in \RR^{2m}$ and $c_i \in \RR^{2ms}$ are the right-hand sides of the Newton-like iteration computed from \eqref{hbvm2}.

Clearly, we shall obtain different linear systems, depending on the boundary conditions in \eqref{ham_bvp}: generally, their efficient solution requires the use of different, specifically tailored, linear solvers. This particular aspect is only sketched here and will be considered elsewhere.

Preliminarly, we observe that simple matrix manipulations would allow us to separate the computation of the stages updates $\Delta_i^{(j)}$ from the solution updates $\delta_i^{(j)}$. In fact, if $K_i$ is nonsingular, from \eqref{struct_matr} one easily derives that $$\Delta_i^{(j)} = K_i^{-1} (b_i - V_i \delta_i^{(j)}), \quad i=1,\dots,n-1.$$ However, due to possible stability problems (indeed, $K_i$ may be ill conditioned or even singular), in general it is preferable to avoid this reduction step, thus solving a linear system whose dimension depends on $m$, $s$, and $n$ (as matter of fact, it turns out to be $\approx\,2m(s+1)(n+1)$).

\subsection{Separated boundary conditions}
The simplest case is when problem \eqref{ham_bvp} is defined by means of $r<2m$ initial and $2m-r$ final nonlinear conditions: $$g_a(y_0)=0\in\RR^r, \qquad g_b(y_n)=0\in\RR^{2m-r}.$$  Then, their linearization provides additional equations in the form
\begin{equation}
\label{sbc}
B_a \delta_0^{(j)} = b_{0a} \in\RR^r, \qquad
B_b \delta_n^{(j)} = b_{0b}\in\RR^{2m-r},
\end{equation}
where $B_a \in R^{r \times 2m}$ and $B_b \in R^{2m-r \times 2m}$. Ordering the two equations as the first and the last one, \eqref{struct_matr} and \eqref{sbc} produce a linear system with an {\em Almost Block Diagonal} (ABD) coefficient matrix \cite{ACFGKRPW}. In such a case, the solution is efficiently obtained by means of direct solvers that generalize the LU factorization (see \cite{ACFGKRPW} for a complete review), with a computational cost consisting into a number of operations proportional to $m^2s^2n$ and no fill-in (i.e., no additional memory is required for the factorization, besides that needed for storing the blocks in the coefficient matrix).

\subsection{Non-separated boundary conditions}
Suppose problem \eqref{ham_bvp} is defined by means of $2m$ (generally nonlinear) boundary conditions involving $y_0$ and $y_n$, i.e., $$g(y_0,y_n)=0\in\RR^{2m}.$$ Then, the linearization of this condition produces the equation
\begin{equation}
\label{nsbc}
B_a \delta_0^{(j)} + B_b \delta_n^{(j)} = b_0\in\RR^{2m},
\end{equation}
that, combined with \eqref{struct_matr}, gives a nonsingular {\em Bordered Almost Block Diagonal} (BABD) linear system whose factorization is conveniently handled by means of a cyclic reduction approach, as is shown in \cite{AR}. This algorithm requires twice the number of operations as in the previous case and generates a fill-in which is essentially equal to $2m(2m+s)n$ memory locations \cite{AR}.

\subsection{Periodic boundary conditions}
\label{sec:3.3}
 From a numerical point of view, the most difficult case to be solved  is when problem \eqref{ham_bvp} is defined with periodic boundary conditions, i.e.,
\begin{equation}
\label{pbc1}
g(y_0,y_n)\equiv y_0 - y_n = 0.
\end{equation}
In such a case, in fact, the continuous problem admits always an infinite number of solutions and, during its discretization, it is necessary to consider additional conditions at $t_0$ (called {\em anchors}) for one or more components of the solution  \cite{NB}, e.g. in the form:
\begin{equation}
\label{pbc2}
B_a y_0 = b_0\in\RR^r,
\end{equation}
where $B_a \in \RR^{r \times 2m}$, with $1\le r\le 2m$. Consequently, the resulting linear system is overdetermined and, hence, it is solved by means of a least square approach. In general, it turns out that, after convergence of the Newton-type iteration, the residual vector is nonzero, though having a $O(h^{2s})$ norm, consistently with the order $2s$ of the underlying HBVM$(k,s)$ method.

For these periodic problems, it could be of interest to regard $t_f$ as unknown\,\footnote{As matter of fact, the period ~$T\equiv t_f-t_0$~ of the orbit could be  not known.} and to locate a periodic orbit by exploiting its energy value instead of its period. Using a constant stepsize $h=\frac{t_f-t_0}n$, this means to retain $h$ as a further unknown and to add an extra equation involving the given value of the Hamiltonian function at the first point, i.e., 
\begin{equation}\label{H0}
H(y_0) = c_0\in\RR.
\end{equation}
To obtain the linear system at the corresponding $j$-th  Newton-like iteration,  we observe that:
\begin{itemize}
\item the unknown $\delta_n^{(j)}$ has been removed  by virtue of (\ref{pbc1}) and, consequently, also the block-row corresponding to \eqref{pbc1} itself;
\item the update $\delta_h^{(j)}$, such that $h_{j+1} = h_j+\delta_h^{(j)}$, has been included, where $h_j$ the the current approximation to the unknown value of the correct stepsize.
\end{itemize} Conseqently, the linear system generated by the Newton-type iteration assumes the form
\begin{equation}
\label{struct_matr2}
\pmatrix {cccccccccc}
B_H \\
B_a \\
V_1 & K_1 &&&&&&& w_1 \\
L_1 & U_1^T & I_{2m} &&&&&& v_1 \\
&& V_2 & K_2 &&&&& w_2 \\
&& L_2 & U_2^T & I_{2m} &&&& v_2 \\
&&&& \ddots & \ddots &&& \vdots \\
&&&&&& V_n & K_n & w_n \\
I_{2m} &&&&&& L_n & U_n^T  & v_n
\endpmatrix
\pmatrix {c}
\delta_0^{(j)} \\ \Delta_0^{(j)} \\ \delta_1^{(j)} \\ \Delta_1^{(j)} \\ \delta_2^{(j)} \\ \vdots \\ \delta_{n-1}^{(j)} \\ \Delta_{n-1}^{(j)} \\ \delta_h^{(j)}
\endpmatrix =
\pmatrix {c}
\gamma \\ \tilde b_0 \\ b_1 \\ c_1 \\ b_2 \\ c_2 \\ \vdots \\  b_n \\ c_n
\endpmatrix,
\end{equation}
where: $B_H = \nabla^T H(y_0^{(j)}) \in \RR^{1 \times s}$; $\gamma\in\RR$ and $\tilde{b}_0\in\RR^r$ derive from the linearization of \eqref{pbc2} and \eqref{H0}, respectively; blocks $V_i,U_i^T,L_i,K_i$ are the same as defined in \eqref{blocks} but with $h_j$ in place of $h$; finally, $w_i$ and $v_i$ are vectors of suitable dimension obtained by differentiating \eqref{hbvm2} with respect to $h$. 

Because of the anchor equation, system \eqref{struct_matr2} still requires a least square approach to be solved: this can be efficiently done by using an algorithm similar to that considered for non-separated boundary conditions \cite{ACFGKRPW,AR}.

\section{Applications to celestial mechanics\\ and astrodynamics}
\label{sec:4}

The following numerical tests are mainly concerned with the motion of a body with negligible mass (planetoid) in the gravitational field generated by two celestial bodies with finite mass (primaries) rotating around their common center of mass in circular orbits. Such a dynamical system is referred to as the \textit{circular restricted three-body problem} (CRTBP) and its interest goes back to the second quarter of the eighteenth century, in the context of the lunar theory. A renewed interest arose starting from the late 1960s up to present day  and is testified by a rich and growing literature on the design and analysis of a variety of orbits connected with the motion of  spacecrafts, satellites and asteroids \cite{Fa,Ho,JM,PGWZ,SKWLMPRW}. 

We here consider the case where the two primaries are the Sun and the Earth+Moon whose masses are denoted $m_1$ and $m_2$. Usually the units are normalized and chosen so that the properties of the resulting dynamical system depend on a single parameter $\mu$, defined as the ratio $\frac{m_2}{m_1+m_2}$. In our situation we have $\mu = 3.04036\cdot10^{-6}$.  To obtain dimensionless coordinates the following normalizing assumptions are introduced:  
\begin{enumerate}
\item the total mass of the system is~  $m_1+m_2=1$; 
\item the unit of length is the distance between the two primaries, i.e., $R=1.49589\cdot10^8$km;
\item the unit of time is $1/n$, where $n=1.99099\cdot10^{-7}$rad/s is the constant angular velocity of the Sun and Earth/Moon around their center of mass $C_M$. 
\end{enumerate}
Notice that, from the above hypotheses, the gravitational constant is unity, $G=1$. It is also common to write down the equations of motion of the planetoid in a frame where the primaries are stationary. This is accomplished by introducing a rotating (synodic) orthogonal frame centered at $C_M$, with the $x$-$y$ axes lying in the plane of the Sun-Earth/Moon orbit, the $x$ axis being oriented from the Sun toward the Earth, and the $z$ axis forming a right-hand frame with the other axes. Thus, the Sun and the Earth are located on the $x$ axis at the abscissae $-\mu$ and $1-\mu$ respectively. 

Let $q(t)=[q_1(t),q_2(t),q_3(t)]^T$ be the coordinates of the planetoid at time $t$ and set $p(t)=[p_1(t),p_2(t),p_3(t)]^T\equiv [\dot q_1(t)-q_2(t),\dot q_2(t)+q_1(t),\dot q_3(t)]^T$ the vector of conjugate momenta. The Hamiltonian function in  non-dimensional form associated with the dynamical system governing the motion of the planetoid is 
\begin{equation}
\label{3ham}
H(q,p)=p_1q_2-p_2q_1+\frac{1}{2}p^Tp-\frac{1-\mu}{r_1}-\frac{\mu}{r_2},
\end{equation}
where $r_1=((q_1+\mu)^2+q_2^2+q_3^2)^{1/2}$ and $r_2=((q_1-(1-\mu))^2+q_2^2+q_3^2)^{1/2}$ are the distances of the planetoid from the Sun and the Earth/Moon respectively. 

It is well-known that such a dynamical system admits five equilibrium points referred to as \textit{Lagrangian} or \textit{libration} points: three ($L_1$, $L_2$, $L_3$) are collinear with
the primaries and the other two ($L_4$ and $L_5$) form an equilateral triangle with them.

Periodic and quasi-periodic orbits around libration points are  suited for a number of mission applications. For example, Sun-Earth libration points are commonly used for deep space or Sun activity observations. In the following experiments we are interested in the dynamics around $L_2$, which is located  beyond the Earth, on the $x$-axis, at the abscissa $1.010075$. 

\medskip
An interesting problem in astrodynamics is the \textit{optimal orbit transfer}, which consists in finding the optimal control laws that drives a spacecraft from an initial state, say $P_1$, to a desired final state $P_2$ in a given time $T$. Here, the term \textit{optimal} means that the amount of propellant needed to produce the change in orbital elements is minimized.  

Since the fuel consumption is proportional to changes in the velocity, an input vector $u(t)=[u_1(t),u_2(t),u_3(t)]^T$ enters the dynamical system to control the acceleration of the vehicle along the three axes. This is accomplished by considering a new non-autonomous Hamiltonian function $$\bar H(q,p) =H(q,p)+q^Tu,$$ where $H(q,p)$ is as in \eqref{3ham}. 
Our optimal control problem is then formulated as follows:
\begin{quote} \em
Minimize the quadratic cost
$
J=\frac{1}{2} \int_0^T ||u(t)||_2^2 \mathrm{d}t,
$
subject to the dynamics induced by $\bar H(q,p)$ and the boundary conditions corresponding to the states $P_1$ and $P_2$.
\end{quote}
We assume that the control input is unconstrained and regular. The Pontryagin maximum principle is often used to attack this problem. Setting $y^T=[q^T,p^T]$ (state variables) and $\lambda=[\lambda_1,\dots,\lambda_6]^T$ (costate variables), one considers the augmented Hamiltonian function $$\widetilde H(y,\lambda,u)=\frac{1}{2}u^Tu+\lambda^TJ\nabla \bar H(q,p).$$ Then, the necessary conditions for optimality are 
$$
\dot y = \frac{\partial \widetilde H}{\partial \lambda}, \qquad  \dot \lambda = - \frac{\partial \widetilde H}{\partial y}, \qquad \frac{\partial \widetilde H}{\partial u}=0.
$$
The third equation gives $u_i=-\lambda_{(3+i)}$, $i=1,2,3$, so that the resulting system is autonomous and only depends on the state and costate variables. It is defined by the Hamiltonian 
\begin{equation}
\label{3ham_extended}
\begin{array}{rl}
\widehat H(y,\lambda)& =\frac{1}{2}(\lambda_4^2+\lambda_5^2+\lambda_6^2)+\lambda^T(J\nabla H(q,p) - [0,0,0,\lambda_4,\lambda_5,\lambda_6]^T)\\[.25cm]
& =\lambda^TJ\nabla H(q,p)-\frac{1}{2}(\lambda_4^2+\lambda_5^2+\lambda_6^2).  
\end{array}
\end{equation}

We now consider a few applications concerning the above problems. All experiments have been carried out in Matlab  (in double precision arithmetic) by using its sparse linear solvers.\footnote{More efficient linear solvers will be studied elsewhere.}

\subsection{Computation of Lyapunov orbits} \textit{Lyapunov orbits} are periodic orbits surrounding a libration point  in the planar CRTBP, where the term \textit{planar} means that the planetoid moves in the same plane as the primaries, namely the $x$-$y$ plane: $q_3(t)=0$, $p_3(t)=\dot q_3(t)=0$. We are interested in the computation of Lyapunov orbits emanating from the point $L_2$, which we here assume as the origin of the axes. Their existence is guaranteed by Lyapunov's center theorem \cite{MHO}, which also states that Lyapunov orbits form a one-parameter family parametrized by the Hamiltonian integral. Thus, it makes sense to search for a Lyapunov orbit by fixing either its period  or its energy level.  We consider both situations and notice that in the latter case an energy-conserving method is more appropriate since it provides a numerical solution that precisely lies on the required energy set.  An analysis of the monodromy matrix associated with Lyapunov orbits shows their instability character, which makes their computation a delicate issue. 

\begin{figure}[t]
\begin{center}
\includegraphics[width=12cm,height=6cm]{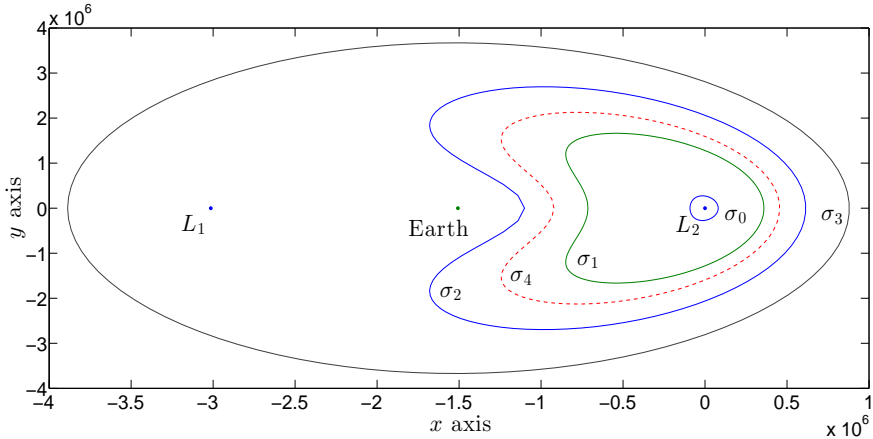} 
\end{center}
\caption{Some Lyapunov orbits surrounding the libration point $L_2$. Their computation may be carried out by passing as input information either their period or their energy level.}
\label{fig00}
\end{figure}

We discretize the time interval into $n=100$ uniform points and use the method HBVM$(6,2)$ which ensures a practical energy conservation for the problem at hand and the used stepsize (see \eqref{errH}). As initial condition for the Newton iteration, we consider a periodic orbit very close to the equilibrium point $L_2$ obtained as the solution of the linearized problem: it is the closed curve labelled as $\sigma_0$ in Figure \ref{fig00} and corresponds to a period $T=178$ days. 

The curve $\sigma_1$ denotes the Lyapunov orbit with period $T=200$ days and has been obtained by considering periodic boundary conditions, as discussed in Section  \ref{sec:3.3}. The energy level associated with this orbit is $H_2\simeq -1.5002604$.

Starting from $\sigma_1$, we attempt to find the Lyapunov orbit corresponding to the energy level $H_3=-1.5001$ and thus we solve the iteration described at \eqref{struct_matr2}. We obtain the orbit labelled as $\sigma_2$ in Figure \ref{fig00}: its  period is $T_3\simeq 251.34$ days. 

The search of $\sigma_2$ via its period $T_3$ rather than its energy level $H_3$ starting from $\sigma_1$ would not provide the desired result: whatever method in the family HBVM$(k,2)$ we choose, the iteration process converges to a different periodic orbit $\sigma_3$ that embraces the Lagrangian points $L_1$ and $L_2$ other than the Earth.  This orbit has period $T_3$ but its energy is $H_3\simeq -1.500177$. To retrieve the correct Lyapunov orbit we need to compute an intermediate curve, such as $\sigma_4$, that corresponds to a period $T_4=220$ days. 

In conclusion, it seems that, for this problem, the continuation technique performs better if based upon the value of the energy rather than of the period.

\subsection{The Hill three-body problem}
The Hill problem is a special, simplified case of the planar CRTBP. It studies the motion of the planetoid in a neighborhood of the Earth, which is conveniently taken as the new origin of the synodic frame via the change of coordinates $q_1 \rightarrow q_1+(1-\mu)$, $q_2\rightarrow q_2$. The assumption on the location of the planetoid permits a simplification of the equations describing its dynamics. Essentially, one discards the terms of order at least three in $q_1$ and $q_2$ in the Taylor expansion of the potential around $(0,0)$, and performs an additional change of variables to simplify the final shape of the equations, making them independent of the parameter $\mu$ (see, for example, \cite{AKN} for details). The Hamiltonian function arising from these transformations reads
\begin{equation}
\label{hill_ham}
H(q_1,q_2,p_1,p_2)=p_1q_2-p_2q_1+\frac{1}{2}(p_1^2+p_2^2)-\frac{1}{(q_1^2+q_2^2)^{1/2}}+\frac{1}{2}q_2^2-q_1^2.
\end{equation}
This reduced system admits only two equilibrium points located on the $x$ axis on both sides of the Earth: $L_1=\left(-(1/3)^{1/3},0\right)$ and $L_2=\left((1/3)^{1/3},0\right)$.

\begin{figure}[t]
\begin{center} 
\includegraphics[width=10cm,height=5cm]{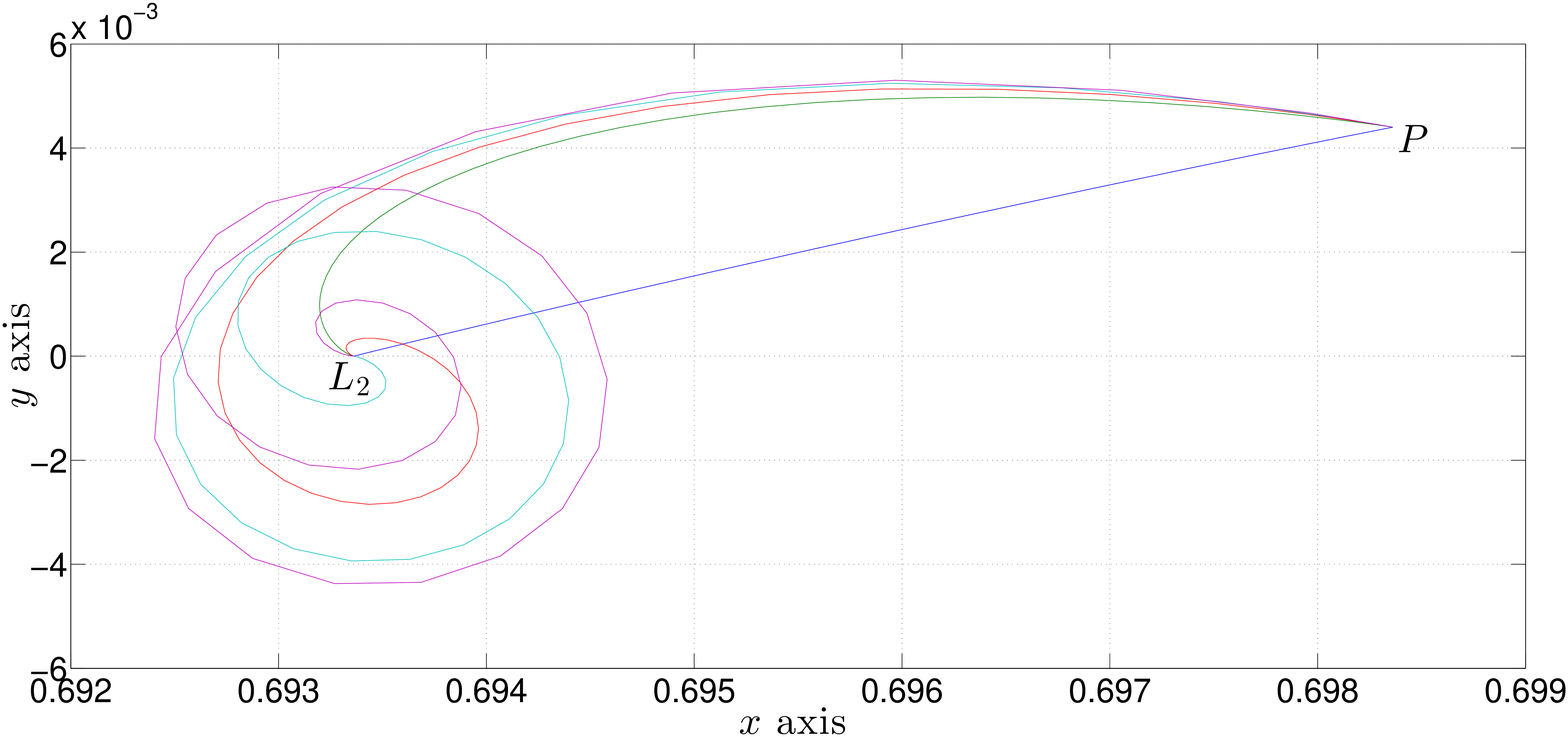} \\
\includegraphics[width=10cm,height=5cm]{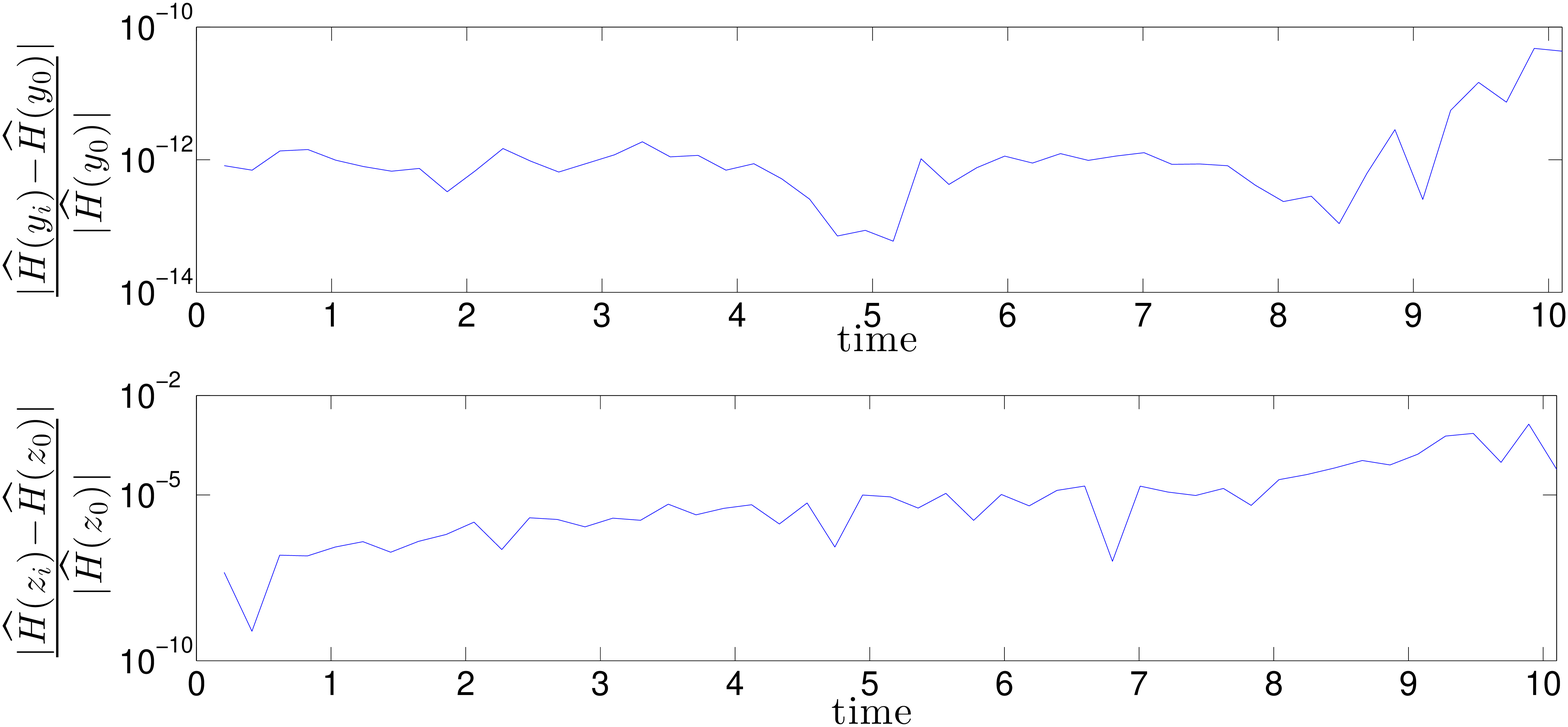}
\end{center}
\caption{Upper picture: orbits of a spacecraft driven from the libration point $L_2$ to a close point $P$ for several transfer times. Lower pictures: relative error in the Hamiltonian function \eqref{3ham_extended} evaluated along the numerical solutions obtained by the HBVM$(4,2)$ (intermediate plot) and the Gauss method of order $4$ (bottom plot). Both solutions correspond to $t_f=8.1$.}
\label{fig0}
\end{figure}

We consider a deployment problem, taken from \cite{GuSc}, consisting in transferring a spacecraft from the point $L_2=((1/3)^{1/3},0)$ to the point $P=((1/3)^{1/3}+0.005, 0.0044)$. In both points the velocity is assumed null and the transfer time is increased as $t_f=0.1,2.1,4.1,6.1,8.1$. Due to the fact that the dynamics takes place near an equilibrium point, we choose the HBVM$(4,2)$ method as integrator, since two silent stages are enough to guarantee a practical energy conservation. The top picture in Figure \ref{fig0} shows the five trajectories of the spacecraft corresponding to the selected transfer times. As $t_f$ is increased, the spacecraft circles around the point $L_2$, in a spiral-shaped orbit, before approaching the final point $P$. The intermediate plot of Figure \ref{fig0} reports the relative error in the Hamiltonian function \eqref{hill_ham} evaluated along the numerical solution $\{y_i\}$ corresponding to $t_f=8.1$. We notice that it is bounded by $10^{-10}$ and cannot be further reduced even if we increase the number of silent stages. This is an effect of the use of finite precision arithmetic and the fact that the order of the Hamiltonian along the orbit is $10^{-6}$. For comparison purposes, in the bottom plot of Figure \ref{fig0}, we have also  included the corresponding error produced by the 4-th order Gauss method (i.e., HBVM(2,2)).        

\subsection{Computation of Halo orbits}
\textit{Halo orbits} are out-of-plane periodic orbits which trace a halo around the Earth.  We are interested in reproducing Halo orbits around the point $L_2$. 

We have implemented the HBVM$(6,2)$ formula and adapted the algorithm in order to compute periodic solutions in the two different situations where we are given either the period ~$T\equiv t_f-t_0$~ of the orbit or its energy level $H_0$. In the latter case, according to what said in Section~\ref{sec:3.3}, the stepsize of integration $h$ is regarded as an extra unknown and the scalar equation $H(q_0,p_0)=H_0$ is added to the set of boundary conditions. 

\begin{figure}[t]
\hspace*{-0.3cm} \includegraphics[width=7cm,height=6cm]{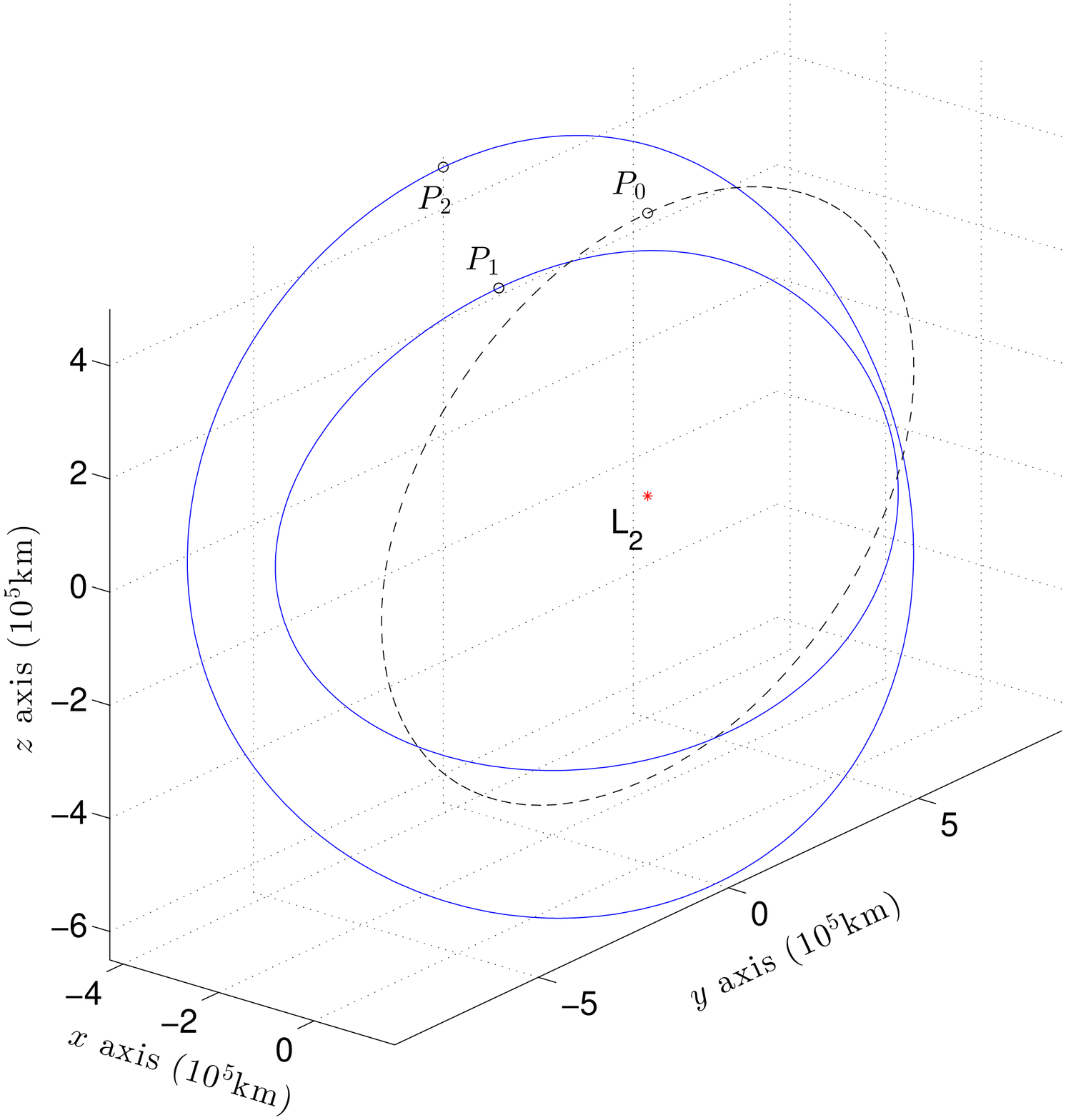} \hspace*{-0.8cm}
\includegraphics[width=6cm,height=6cm]{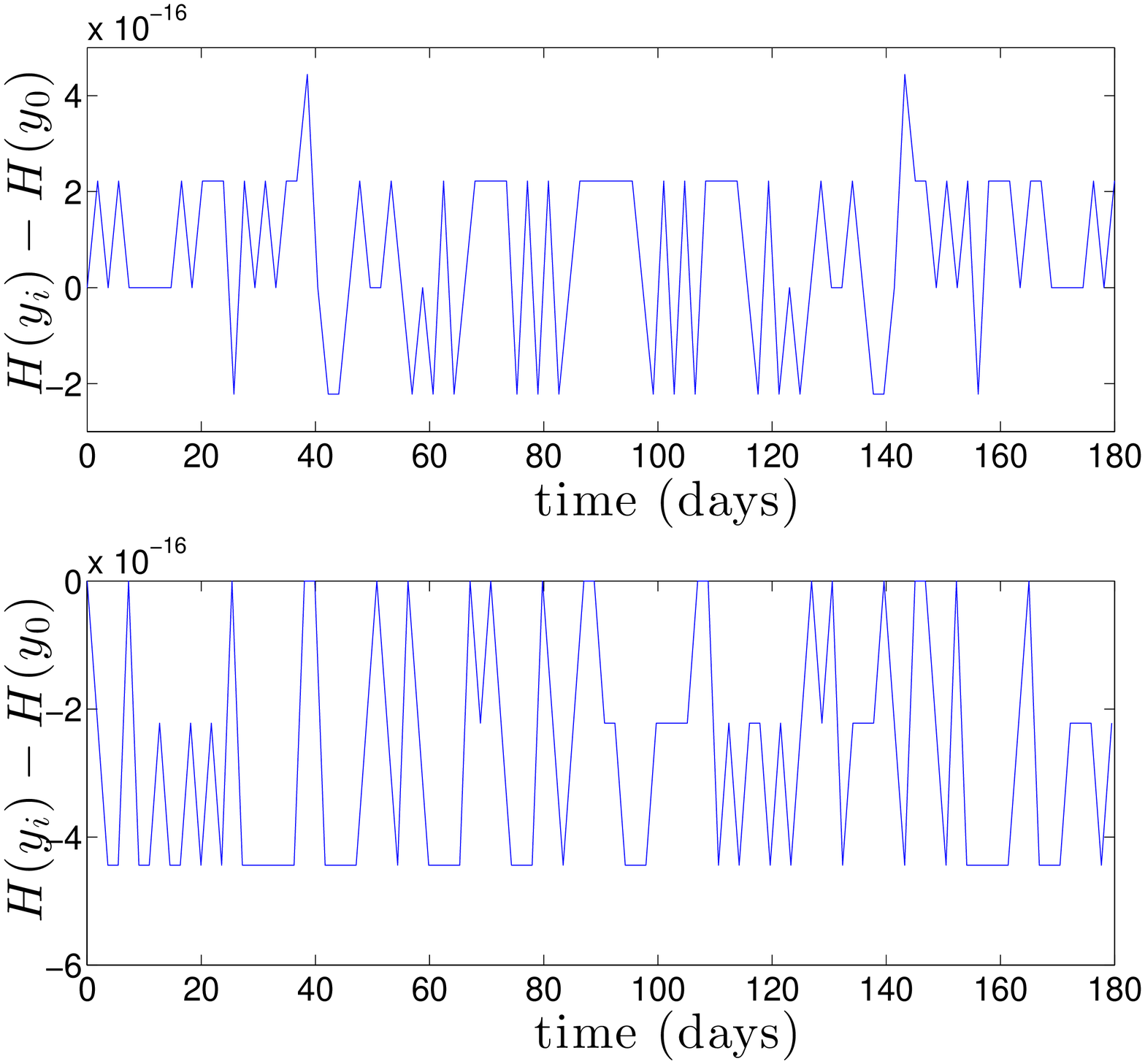}
\caption{Left picture: two halo orbits around the libration point $L_2$ (solid lines) and the initial guess for the Newton iteration scheme associated to the method (dashed line). Right picture: the Hamiltonian function \eqref{3ham} is precisely conserved along the numerical solutions.}
\label{fig1}
\end{figure}

In both cases, an elliptic curve lying on a plane orthogonal to the $x$ axis and passing through $L_2$ has been chosen as initial guess for the Newton iteration. The starting (and ending) point $P_0$ of this curve has been set at the upper end of the vertical axis of the ellipse (see Figure \ref{fig1}). The number of points in the numerical approximation is $n=100$. 

The left picture of Figure~\ref{fig1} displays the initial guess (dashed line) together with two halo orbits (solid lines). The inner one is the halo orbit corresponding to a period $T_1=180$ days. The energy level of this first numerical approximation is $H_1\approx -1.500394$. Conversely, the outer halo orbit has been computed on the basis of its energy level, which has been set to $H_2=-1.50036$. Notice that in the non-dimensional system $H_2\approx H_1(1+2\cdot 10^{-5})$ while the actual distance of the topmost points of the two orbits is $\overline{P_1P_2}=2\cdot 10^5$km. The period corresponding to the energy level $H_2$ is $T_2=179.19$ days. 
The right pictures of Figure \ref{fig1} show that the energy error is close to the machine precision in both cases.

\begin{figure}[t]
\hspace*{-0.3cm} \includegraphics[width=7cm,height=6cm]{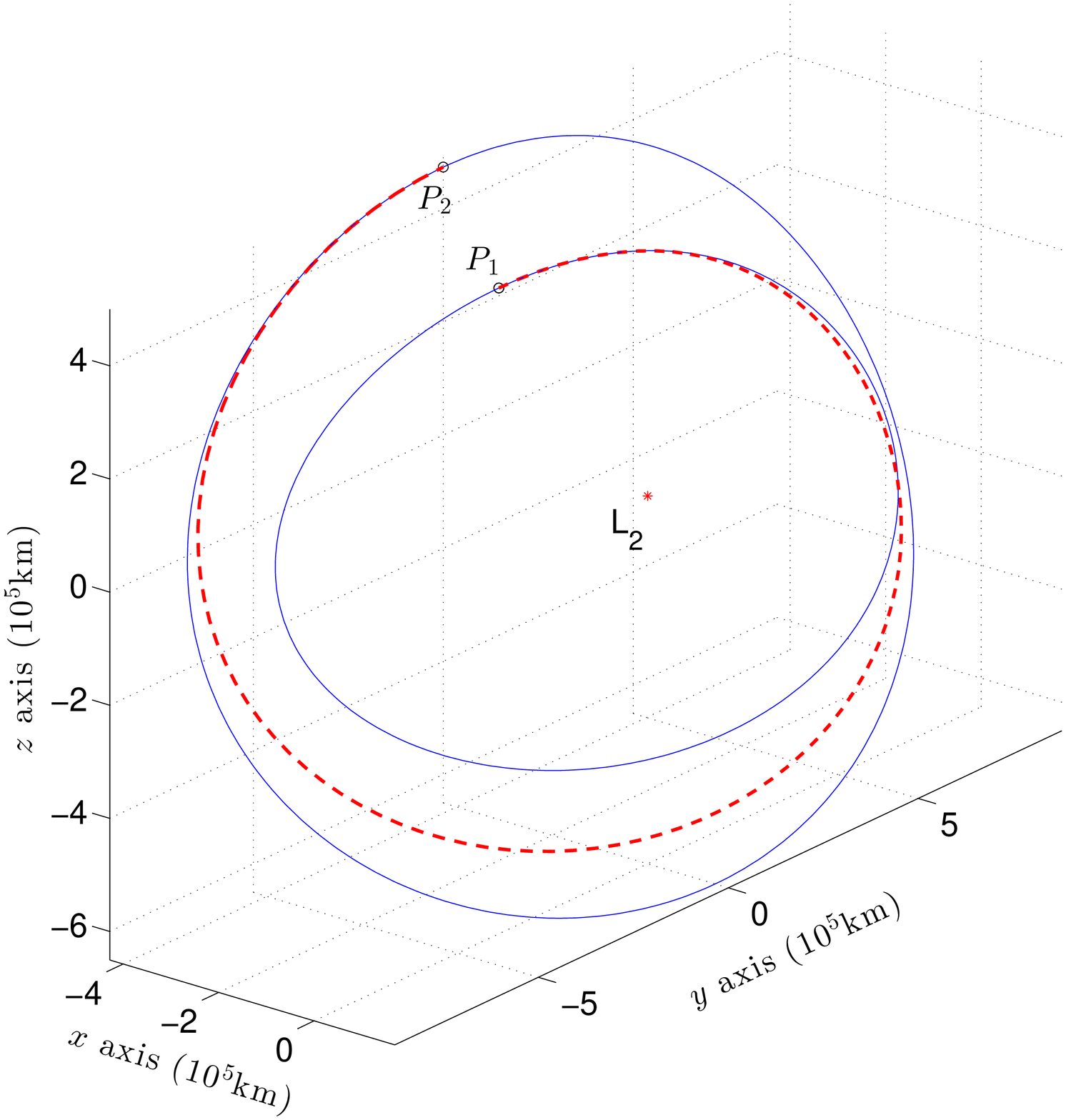} \hspace*{-0.8cm}
\includegraphics[width=6cm,height=6cm]{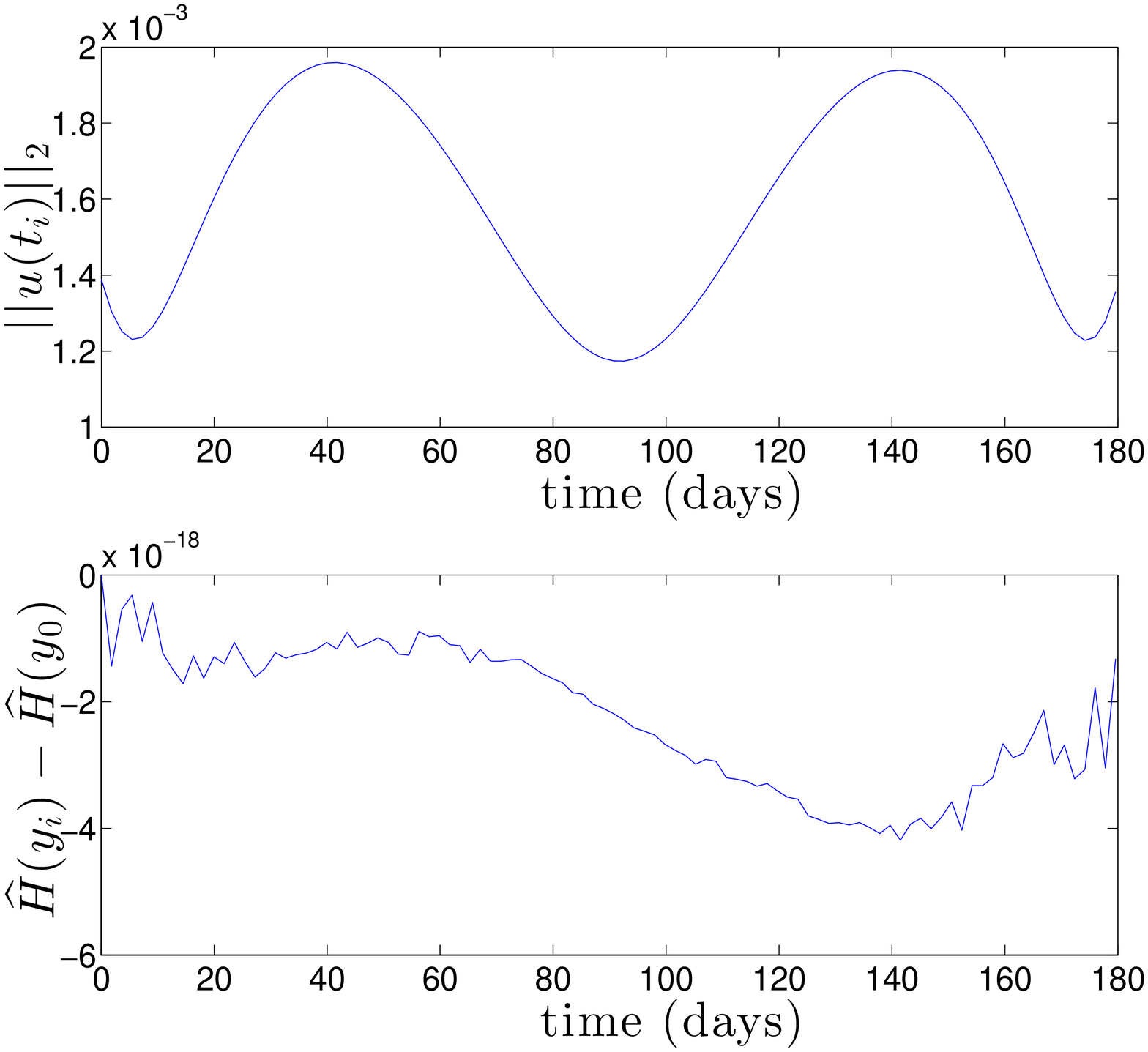}
\caption{Left picture: optimal orbit transfer between  two halo orbits (dashed line). Right picture: norm of the optimal control variable $u(t)$ (upper plot) and error in the Hamiltonian function \eqref{3ham_extended} (lower plot).}
\label{fig2}
\end{figure}

We also consider the optimal orbit transfer problem consisting in transferring a spacecraft from the inner to the outer halo orbit and, specifically, from the point $P_1$ to the point $P_2$. In Figure \ref{fig2} we show the optimal control trajectory joining the points $P_1$ and $P_2$ in a time $T=(T_1+T_2)/2$ (left picture, dashed line) together with the norm of the optimal control variable $u(t)$ and the error $\widehat H(y_n)-\widehat H(y_0)$ in the Hamiltonian \eqref{3ham_extended} (right picture). 

\section{Conclusions}
\label{sec:5}
In this paper, we have extended the use of HBVMs to the solution of Hamiltonian Boundary Value Problems. HBVMs form a subclass of Runge-Kutta methods,  characterized by a rank-deficient coefficient matrix, that provide a numerical solution along which the Hamiltonian function is precisely conserved. Their implementation has been adapted in order to handle different kinds of boundary conditions. In particular, separate and periodic boundary conditions arise in several problems of celestial mechanics and astrodynamics, such as the periodic orbit detection and the optimal spacecraft orbit transfer. A few numerical tests in this direction have shown the good potentialities of the methods.




\end{document}